\documentclass[10pt,twoside]{amsart}

\usepackage{enumerate}
\usepackage[all]{xy}
\usepackage[hyperindex=true,plainpages=false,colorlinks=false,pdfpagelabels]{hyperref}

\theoremstyle{plain}
\newtheorem{The}{Theorem}
\newtheorem*{The*}{Theorem}
\newtheorem{Pro}[The]{Proposition}
\newtheorem{Lem}[The]{Lemma}

\newtheorem*{Cor*}{Corollary}

\theoremstyle{definition}
\newtheorem*{Def}{Definition}

\theoremstyle{remark} 
\newtheorem{Rem}[The]{Remark}

\newtheorem*{Rem*}{Remark}

\numberwithin{equation}{section}

\DeclareMathOperator{\End}{End}
\DeclareMathOperator{\Hom}{Hom}

\DeclareMathOperator{\SL}{SL}

\DeclareMathOperator{\SU}{SU}
\DeclareMathOperator{\Tr}{tr}             

\DeclareMathOperator{\Id}{Id}

\DeclareMathOperator{\vol}{vol}

\renewcommand{\Im}{\operatorname{Im}}

\newcommand{\dvector}[1]{{\left(\begin{matrix}#1\end{matrix}\right)}}
\DeclareMathOperator{\dbar}{\bar\partial}
\DeclareMathOperator{\del}{\partial}

\newcommand{\R}{\mathbb{R}}
\newcommand{\C}{\mathbb{C}}

\newcommand{\Z}{\mathbb{Z}}

\renewcommand{\H}{\mathbb{H}}  

\newcommand{\CP}{\mathbb{CP}}

\begin{document}

\title[Lawson's genus two minimal surface and meromorphic connections]{Lawson's genus two minimal surface and meromorphic connections}

\author{Sebastian Heller}


\address{Sebastian Heller \\Institut f\"ur Mathematik\\  Universit{\"a}t T\"ubingen\\ Auf der Morgenstelle10\\ 72076 T¬ubingen\\ Germany }

\email{heller@mathematik.uni-tuebingen.de}



\date{\today}

\thanks{Authors supported by SFB/Transregio 71}

\begin{abstract} 
We investigate the Lawson genus $2$ surface by methods from integrable system theory. We prove that the associated family of flat connections comes
from a family of flat connections on a $4-$punctured sphere. We describe the symmetries of the holonomy and show that it is already determined 
by the holonomy around one of the punctures.
We show the existence of a meromorphic DPW potential for the Lawson surface 
which is globally defined on the surface. We determine this potential explicitly up to two unknown functions depending only on the spectral parameter.
 \end{abstract}

\maketitle


\section{Introduction}
\label{sec:intro}
Concrete examples of the right type have always been fruitful in mathematics.
The construction of constant mean curvature (CMC) tori by Wente \cite{W} has stimulated the work on CMC tori in $3-$dimensional space forms by many authors. After Abresch's \cite{A} analytical description of the Wente tori, a complete classification of CMC tori in terms of holomorphic data was given by Pinkall and Sterling \cite{PS} and Hitchin \cite{H1} independently. This lead to the construction of all CMC 
tori in terms of theta functions by Bobenko \cite{B2}.

There are also examples of compact minimal surfaces and CMC surfaces in $\R^3$ of higher genus, see \cite{L}, \cite{KPS} or \cite{K}.
The genus $2$ minimal surface $M\subset S^3$ of Lawson \cite{L}, which we are going to study here, might be the most simple one. But none of these surface is known explicitly and the construction of them gives no hint how to describe all compact CMC surfaces in space forms. The aim of this paper 
is to study Lawson's genus $2$ minimal surface $M$ in a more explicit way. The hope is, that this provides some insight into a theory of higher genus surfaces. 

There is
a general method due to Dorfmeister, Pedit and Wu, \cite{DPW}, which produces, in principal, all CMC surfaces ( and, more generally, harmonic maps into symmetric spaces). A CMC surface in a $3-$dimensional space form
can be described by their associated family of flat connections $\nabla^\zeta$ on a complex rank $2$ bundle $V.$
The idea of the DPW method is to gauge $\nabla^\zeta$ into a family of meromorphic connections of a special form, the so-called DPW potential, in a way which can be
reversed. The advantage is that one can write down meromorphic connections easily.
On simply connected domains, each minimal surface can be obtained from such a family of meromorphic connections. To obtain a surface
one takes a $\zeta-$depending parallel hence holomorphic frame and splits it into the
unitary and positive parts via Iwasawa decomposition in the loop group. Then the unitary part is a parallel frame of a family of unitary connections 
describing a minimal surface. The surface obtained in this fashion depends on the $\zeta-$depending initial condition of the parallel frame.
Dressing, i.e. changing this initial condition, will give different surfaces.
If one wants to make surfaces with topology via DPW, one has to
ensure that one can patch simply connected domains together. 
This has been worked out only in very special cases, for example for trinoids, the genus zero CMC surfaces with three Delauney ends, or CMC tori. 
Up to now there are no examples of closed higher genus surfaces. We show how the DPW method can be applied to the case of the Lawson surface $M,$ and prove that a globally defined DPW potential
for the Lawson surface does exists on $M.$ We determine this potential almost explicitly.

In the first part we recall the gauge theoretic description of minimal surfaces in $S^3.$ We give an explicit link to the local description of minimal surfaces via the extended frame. In the third section we shortly explain Lawson's construction of compact minimal surfaces in $S^3.$ We collect all the symmetries and all holomorphic data of Lawson's genus $2$ surface. Especially, we will determine the spin bundle $S,$ and we show that the associated rank $2$ bundle $V$ with the holomorphic structure $(\nabla^0)''$ is stable.

The fourth section is devoted to the study of the holonomy of the Lawson surface $M.$ By using the $\Z_3-$symmetry of the Lawson surface we show that the associated family of flat connections $\nabla^\zeta$ comes naturally from a family of flat connections globally defined on $\mathcal O(1)\oplus\mathcal O(-1)\to\CP^1$ with singularities at the $4$ branch points
of $\pi\colon M\to M/\Z_3=\CP^1.$ The holonomy of the Lawson surface is therefore given by the holonomy of a family of flat connections on a $4-$punctured sphere (see theorem \ref{holonomyupdown}). Using the other symmetries we prove that this holonomy  is entirely given by the monodromy of the family of connections
on $\CP^1$ based at $0\in\CP^1$ around any of its four singularities.

The last part deals with the construction of a DPW potential for the Lawson surface. We prove that one can find a globally defined gauge with pole like singularities at the Weierstrass points of $M,$ such that the family of connections obtained by gauging is a meromorphic family of connections with respect to the fixed direct sum spin holomorphic structure $S^*\oplus S$ on $V.$ The gauge is positive
in the loop group, i.e. it extends to $\zeta=0$ in a special form, so that one
can get back the Lawson surface by the DPW method. Using the symmetries of the surface, we can show that the DPW potential has corresponding symmetries. 
In fact, there exists a corresponding family of meromorphic connections on $\mathcal O(1)\oplus\mathcal O(-1)\to\CP^1$ with regular singularities at the branch points, and apparent singularities at the images of the Weierstrass points. Moreover,
the symmetries are enough to determine the DPW potential on $M$ (and on $\CP^1$) up to two unknown functions, the accessory parameters, depending on $\zeta$ only (see theorem \ref{DPWex}). These two functions are almost determined by the properties that the holonomy is unitary and that the resulting surface has all symmetries.

The author thanks Aaron Gerding, Franz Pedit and Nick Schmitt for helpful
discussions.
 
\section{Minimal Surfaces in $S^3$}
We shortly describe a way of treating minimal surfaces in $S^3$ 
due to Hitchin \cite{H1}. For more details, one can also consult \cite{He}. 

We consider the round $3-$sphere $S^3$ with its tangent bundle trivialized
by left translation 
$TS^3=S^3\times \Im\H$ and Levi Civita connection given, with respect to the above trivialization, by $\nabla=d+\frac{1}{2}\omega.$
Here $\omega\in\Omega^1(S^3,\Im\H)$ is the Maurer-Cartan form of
$S^3$ which acts via adjoint representation. 
It is well-known that
$S^3$ has a unique spin structure.
We consider the associated complex
spin bundle $$V=S^3\times\H$$
with complex structure given by right multiplication with $i\in\H.$
There is a complex hermitian metric $(.,.)$ on it given by the trivialization and by the identification
 $\H=\C^2.$
The Clifford multiplication is given by
$$TS^3\times V\to V; (\lambda, v)\mapsto \lambda v$$
where $\lambda\in\Im \H$ and $v\in\H.$ This is clearly complex linear.
The induced complex unitary connection is given by
\begin{equation}\label{nablaS3}
\nabla=\nabla^{spin}=d+\frac{1}{2}\omega,
\end{equation}
where the $\Im\H-$valued Maurer-Cartan form acts by left multiplication in the quaternions. Via this construction the tangent bundle
$TS^3$ identifies as the skew symmetric 
trace-free complex linear endomorphisms of $V.$
 
Let $M$ be a Riemann surface and $f\colon M\to S^3$ be a conformal 
immersion. Then the pullback $\phi=f^*\omega$ of the Maurer-Cartan form satisfies the structural equations
$$d\phi+\frac{1}{2}[\phi\wedge\phi]=0.$$
Another way to write this equation is
\begin{equation}\label{structure}
d^\nabla\phi=0,
\end{equation}
where $\nabla=f^*\nabla=d+\frac{1}{2}\phi,$ with $\phi\in\Omega^1(M;\Im\H)$ acting via adjoint representation. 
From now on we only consider the case of $f$ being minimal. Under the assumption of $f$ being
conformal $f$ is minimal if and only if it is harmonic. This is
exactly the case when
\begin{equation}\label{harmonic}
d^\nabla*\phi=0.
\end{equation}
The complex rank $2$ bundle $V:=f^*V\to M$ can be used to
rewrite the equations:
Consider $\phi\in\Omega^1(M;f^*TS^3)\subset
\Omega^1(M;\End_0(V))$ via the interpretation of
$TS^3$ as the bundle of trace-free skew hermitian
endomorphisms of $V.$ Then
$$\frac{1}{2}\phi=\Phi-\Phi^*$$
decomposes
into $K$ and $\bar K$ parts, i.e. 
$\Phi=\frac{1}{2}(\phi-i*\phi)\in\Gamma(K\End_0(V))$ 
and $\Phi^*=\frac{1}{2}(\phi+i*\phi)\in\Gamma(\bar K\End_0(V)).$
Moreover, $f$ is conformal if and only if $\Tr\Phi^2=0.$
In view of a rank $2$ bundle $V$ and $\Tr\Phi=0$ this is equivalent to
\begin{equation}\label{detPhi}
\det\Phi=0.
\end{equation}
Note that $f$ is an immersion if and only if $\Phi$ is nowhere vanishing. 
The equations
\ref{structure} and \ref{harmonic} are equivalent to
\begin{equation}\label{dbarPhi}
\nabla''\Phi=0,
\end{equation}
where $\nabla''=\frac{1}{2}(d^\nabla+i*d^\nabla)$ 
is the underlying holomorphic structure of the pull-back of the spin connection on $V.$
Of course equation \ref{dbarPhi} does not contain
the property that $\nabla-\frac{1}{2}\phi=d$ is trivial. 
Locally, or on simply 
connected sets, this is equivalent to
\begin{equation}\label{curvature}
F^\nabla=[\Phi\wedge\Phi^*]
\end{equation}
as one easily computes. 

Conversely, given an unitary rank $2$ bundle $V\to M$ over
a simply connected Riemann surface with a special unitary connection
$\nabla$ and a trace free field $\Phi\in\Gamma(K\End_0(V))$ without 
zeros, which satisfy \ref{detPhi}, \ref{dbarPhi} and \ref{curvature}, we 
get a conformally immersed  minimal surface as follows:
By equation \ref{dbarPhi} and \ref{curvature}, the unitary connections
$\nabla^L=\nabla-\Phi+\Phi^*$ and $\nabla^R=\nabla+\Phi-\Phi^*$
are flat. Because $M$ is simply connected they are gauge equivalent.
Due to the fact that $\Tr\Phi=0,$ the determinant bundle $\Lambda^2V$ is trivial with respect to all these connections. Hence, the gauge 
is $SU(2)=S^3-$valued with differential $\phi=2\Phi-2\Phi^*.$ Thus it is 
a conformal immersion. The harmonicity follows from \ref{dbarPhi}.
  
From equation \ref{dbarPhi} and \ref{curvature} one sees that
the associated family of connections
\begin{equation}\label{nablafamily}
\nabla^\zeta:=\nabla+\zeta^{-1}\Phi-\zeta\Phi^*
\end{equation}
is flat for all $\zeta\in\C^*.$ As we have see, this family contains all the informations about the surface. It is 
often much easier to describe the family of connections than the minimal surface explicitly, for example in the case of tori, see \cite{H1}
or \cite{B2}, or in the case of a $3-$punctured sphere, see \cite{KKRS}.
The aim of this paper is to study the associated family of flat connections for the Lawson genus $2$ surface, which will be done in section \ref{sec:holonomy} and \ref{DPWpotential}.
  
The geometric significance
of the spin structure of an immersion
$f\colon M\to S^3$ is described in Pinkall
\cite{P}.
We consider the bundle $V$ with its
holomorphic structure $\dbar:=\nabla''.$
As we have seen the complex part $\Phi$ of the differential of a 
conformal minimal surface satisfies $\Tr\Phi=0$ and $\det\Phi=0,$
but is nowhere vanishing.  We obtain a well-defined
holomorphic line subbundle
$$L:=\ker\Phi\subset V.$$
Because $\Phi$ is nilpotent the image of $\Phi$ 
satisfies $\Im\Phi\subset K\otimes L.$ Consider 
the holomorphic section
$$\Phi\in H^0(M;\Hom(V/L,KL))$$
without zeros. The holomorphic structure $\dbar-\Phi^*$ turns $V\to M$ 
into the holomorphically trivial bundle $\underline\C^2\to M.$ As $\Tr\Phi^*=0,$ the determinant line bundle $\Lambda^2V$ of $(V,\dbar)$ is holomorphically trivial. This implies $V/L=L^{-1}$ and we obtain
$$\Hom(V/L,KL)=L^2K$$ as holomorphic line bundles. Because
$L^2K$ has a holomorphic section $\Phi$ without zeros, we get
$$L^2=K^{-1}.$$
Hence, its dual bundle $S=L^{-1}$ is a spinor bundle of the 
Riemann surface $M.$ Clearly, $S^{-1}$ is the only $\Phi-$invariant
line subbundle of $V.$ Moreover, one can show that  $S^{-1}$ is the
$-i-$eigenbundle of the complex quaternionic structure $\mathcal J$ given by quaternionic right multiplication with the right normal vector 
$R\colon M\to S^2\subset\Im\H,$ see \cite{BFLPP} and \cite{He}.
This shows that $S$ gives the spin structure of the immersion.

Let $V= S^{-1}\oplus S$ be the unitary decomposition.
With respect to this decomposition the pull-back of the spin connection on $S^3$ can be written as
\begin{equation}\label{connection_rank2}
\nabla=\dvector{&\nabla^{spin^*} & -\frac{i}{2} Q^*\\ & -\frac{i}{2} Q & \nabla^{spin}},
\end{equation}
where $\nabla^{spin}$ is the spin connection corresponding to the Levi-Civita connection on $M,$  $Q\in H^0(M,K^2)$ is the Hopf field of the immersion $f,$
and $Q^*\in\Gamma(M,\bar K K^{-1})$ is its adjoint, see \cite{He} for details.

The Higgsfield $\Phi\in H^0(M,K\End_0(V))$ can be identified with $$\Phi=1\in H^0(M;K\Hom(S,S^{-1})),$$ and its adjoint $\Phi^*$ is given by the volume form $\vol$ of the induced Riemannian metric.

Let $U\subset M$ be a simply connected open subset and $z\colon U\to\C$ be a holomorphic chart. Write $g=e^{2u}|dz|^2$ for a function $u\colon U\to\R.$
Choose a local holomorphic section $s\in H^0(U;S)$ with $s^2=dz,$ and let $t\in H^0(U,S^{-1})$ be its dual holomorphic section. 
Then $$(e^{-u/2}t,\ e^{u/2} s)$$ is a special unitary frame of $V=S^{-1}\oplus S$ over $U.$ Write the Hopf field $Q=q (dz)^2$ for some local holomorphic function $q\colon U\to\C.$ 

The connection form of the spin connection $\nabla^{spin}$ on the spin bundle $S\to M$ with respect to the local frame $s$ is given by
$-\del u,$ and with respect to $e^{u/2} s,$ it is given by $\frac{1}{2}i*du.$ From formula \ref{connection_rank2} the connection form of $\nabla$ with respect to 
$( e^{-u/2}t,\ e^{u/2} s)$ is
$$\dvector{&-\frac{1}{2}i*du & -\frac{i}{2} e^{-u}\bar q d\bar z\\ & -\frac{i}{2} e^{-u}q d z& \frac{1}{2}i*du}.$$
The Higgsfield $\Phi$ and its adjoint $\Phi^*$ are given by 
$$\Phi=\dvector{&0 &  e^{u}d z\\ & 0& 0},\,\,\,\,\,\, \Phi^*=\dvector{&0 & 0\\ &  e^{u}d \bar z& 0}$$
with respect to the frame $( e^{-u/2}t,\ e^{u/2} s).$
These formulas are well-known, see \cite{DH}, or, in slightly different notation, \cite{B1}.

\section{Lawson's genus $2$ surface}
We recall Lawson's construction \cite{L} of the genus $2$ minimal surface $f\colon M\to S^3.$
We describe the symmetries of the this surface. We use these symmetries to determine the underlying Riemann surface and the holomorphic structures on its associated bundle.
Most of this is well-known, but our arguments in the next sections rely on this.

\subsection{Construction of the Lawson surface}
For two points $A,B\in S^3$ with distance 
$dist(A,B)<\pi$ we denote by $AB$ the minimal oriented geodesic from $A$ to $B.$ If $A$ and $B$ are antipodal, i.e. $dist(A,B)=\pi,$ and $C\in S^3\setminus\{A,B\},$ we denote by $ACB$ the unique oriented minimal geodesic from $A$ to $B$ 
through the point $C.$ 
For a geodesic $\gamma$ and a totally geodesic sphere $S$ we denote the reflections across $\gamma$ and $S$ by $r_\gamma$ and $r_S,$ respectively. 

Let $M$ be an oriented surface with boundary $\gamma,$ and
complex structure $J.$ Let $\gamma$ be oriented and $X\in T_p\gamma$ with $X>0.$ We say 
$\gamma$ represents the oriented boundary if $JX\in T_pM$ represents the exterior normal of the surface for all 
$p\in\gamma\subset M.$

Consider the round $3-$sphere
$$S^3=\{(z,w)\in\C^2\mid |z|^2+|w|^2=1\}\subset\C\oplus\C$$
and the geodesic circles
$C_1=S^3\cap ( \C\oplus\{0\})$ and $C_2=S^3\cap (\{0\}\oplus\C)$ on it.
Take the six points $$Q_k=(e^{i\frac{\pi}{3}(k-1)},0)\in C_1$$ in equidistance on $C_1,$ and the four points $$P_k=(0,e^{i\frac{\pi}{2}(k-1)})\in C_2$$ in equidistance on $C_2.$
Consider the closed geodesic convex polygon $\Gamma=P_1Q_2P_2Q_1$ in $S^3$ with vertices $P_1,Q_1,P_2,Q_2$ and oriented 
edges $P_1Q_2,$ $Q_2P_2,$ $P_2Q_1,$ and $Q_1P_1.$ Then there exists an unique solution for the Plateau problem with boundary 
$\Gamma,$ i.e. a smooth surface which is area 
minimizing under all surfaces with boundary $\Gamma.$
This surface is the fundamental piece of the Lawson surface.
One can reflect this solution at the geodesic $P_1Q_1$ to obtain a smooth surface with piecewise smooth boundary given by the polygon $P_1Q_6P_4Q_1P_2Q_2P_1.$ The  surface obtained in this way can be rotated around $P_1P_2$ by $\frac{2}{3}\pi$ two times, to obtain a new minimal surface, call it $R,$ with possible 
singularity at $P_1$, and with oriented boundary given by the 
oriented edges
$P_2Q_1P_4,$ $P_4Q_6P_2,$ $P_2Q_5P_4,$ $P_4Q_4P_2,$
$P_2Q_3P_4,$ and $P_4Q_2P_2.$
As Lawson has proven, the point P$_1$ is a smooth point on this surface. Now, one can continue, and reflect the resulting surface across the geodesic $C_1.$ Again, the surface 
$R\cup r_{C_1}(R)$ obtained in this way is smooth at each of its points. Moreover it is embedded and orientable. 
The surface is closed as one can see as follows:
The $Q_k$ are fixpoints of $r_{C_1},$ and $r_{C_1}$ interchanges
$P_1$ and $P_3,$ $\ P_2$ and $P_4.$ Moreover this reflection
acts orientation preserving on the surface.
Therefore the oriented boundary edges
$P_2Q_1P_4,$ $P_4Q_6P_2,$ ..., 
$P_4Q_2P_2$ of $R$ are mapped to the oriented
boundary edges 
$P_4Q_1P_2,$ $P_2Q_6P_4,$ ...,
$P_2Q_2P_4$ of $r_{C_2}(R).$ But by the meaning 
of the boundary orientation described above
one sees that $R\cap r_{C_1}(R)$ is closed.

It is proven by Lawson that the zeros
of the Hopf differential $Q$ are exactly at the points $P_1,..,P_4$ of order $1.$ 
Hence, the genus of the Lawson surface is $2$ by Riemann-Roch.

\subsection{Symmetries of the Lawson surface}\label{sym}
There are two types of symmetries of the Lawson surface:
The first type consists of the symmetries (i.e. reflections at geodesics) which were used to construct
the Lawson surface from the fundamental piece. It is clear that they give rise to isometries of the surface.
 The other symmetries are isometries of $S^3$ which map the polygon 
$\Gamma$ to itself. Then, by the uniqueness of the Plateau solution, they give rise to isometries of the Lawson surface, too.

 A generating system 
of the symmetry group of the Lawson surface is given by
\begin{itemize}
\item the $\Z^2-$action generated by $\Phi_2$ with $(a,b)\mapsto(a,-b);$ it is orientation preserving on the surface and its fix points are $Q_1,..Q_6;$
\item the $\Z_3-$action  generated by the rotation $\Phi_3$ around $P_1P_2$ by $\frac{2}{3}\pi,$ i.e.
$(a,b)\mapsto(e^{i\frac{2}{3}\pi}a,b),$ 
which is holomorphic on $M$ with fix points $P_1,..,P_4;$
\item the reflection at $P_1Q_1,$ which is 
antiholomorphic; it is given by $\gamma_{P_1Q_1}(a,b)=(\bar a,\bar b);$
\item the reflection at the sphere $S_1$ corresponding to the 
real hyperplane spanned by $(0,1), (0,i),(e^{\frac{1}{6}\pi i},0),$ with $\gamma_{S_1}(a,b)=(e^{\frac{\pi}{3}i}\bar a, b);$
it is antiholomorphic on the surface,
\item the reflection at the sphere $S_2$ corresponding to the 
real hyperplane spanned by $(1,0), (i,0),(0,e^{\frac{1}{4}\pi i}),$ which is antiholomorphic on the surface and satisfies
$\gamma_{S_2}(a,b)=(a,i\bar b);$
\end{itemize}
Note that all these actions commute with the 
$\Z_2-$action. The last two fix the polygon 
$\Gamma.$ They and the first two map the oriented normal to itself. The third one maps the oriented normal to its negative.

\subsection{The Riemann Surface}\label{RS}
Using the symmetries, one can determine the Riemann surface structure of 
the Lawson surface $f\colon M\to S^3.$ One way to describe the Riemann surface structure is to factor out the $\Z_2-$action
which is exactly the hyperelliptic involution of the genus $2$ surface. Instead of doing this we factor out the $\Z_3-$symmetry
which will be much more useful later on.

The quotient $M/\Z_3$ has an unique structure of a Riemann surface such that $\pi\colon M\to M/\Z_3$ is holomorphic.
The degree of this map is $3$ and its fixpoints are $P_1,..,P_4$ with branch order $2.$ Thus $M/\Z_3=\CP^1$ by Riemann-Hurwitz.
We fix this map by the properties $\pi(Q_1)=0,$ $\pi(P_1)=1$ and $\pi(Q_2)=\infty\in\CP^1.$ Then we have $\pi(Q_3)=\pi(Q_5)= 0,$ and
$\pi(Q_4)=\pi(Q_6)= \infty$ automatically. A symmetry $\tau$ on $M$ gives rise to an action on $\CP^1=M/\Z_3$ if and only if
$\tau(p)$ and $p$ lie in the same $\Z_3-$orbit for all $p\in M.$ This happens for all symmetries described above.

The symmetry $\Phi_2$ defines a holomorphic map $\Phi_2\colon\CP^1\to\CP^1$ which fixes $0$ and $\infty$ and satisfies $\Phi^2=\Id,$
thus $\Phi_2(z)=-z.$ In particular we have $\pi(P_3)=-1.$ Similarly, the induced action of $\gamma_{S_2}$ is antiholomorphic on $\CP^1,$
fixes $0$ and $\infty$ and satisfies $\gamma_{S_2}^2=\Phi_2.$ Therefore $\gamma_{S_2}(z)=\pm i\bar z.$ In fact
$\gamma_{S_2}(z)= i\bar z,$ and we obtain $\pi(P_2)=i$ and $\pi(P_4)=-i.$ 

We collect the symmetries induces on $\CP^1:$
\begin{itemize}
\item the $\Z_2-$action induces $z\mapsto-z;$
\item the reflection at $P_1Q_1$ induces the antiholomorphic map $z\mapsto \bar z;$
\item the reflection at the sphere $S_1$ gives $z\mapsto \frac{1}{\bar z};$ 
\item the reflection at the sphere $S_2$ gives rise to the antiholomorphic map 
$z\mapsto i\bar z.$
\end{itemize}

These observations easily imply the first part of
\begin{Pro}\label{LawsonRS}
The Riemann surface $M$ underlying the Lawson genus $2$ surface is the compactification of the Riemann
surface given by
$$y^3=z^4-1.$$
The Hopf differential of the Lawson genus $2$ surface is given by
$$Q=\pi^*\frac{ir}{z^4-1}(dz)^2$$
for a nonzero real constant $r\in\R.$ \end{Pro}
\begin{proof}
The Hopf differential is $\Z_3-$invariant and has simple zeros at
 $P_1,..,P_4.$ 
Therefore $Q$ is a non-zero complex multiple of 
$\pi^*\frac{1}{z^4-1}(dz)^2.$ 
The Hopf differential
is the $K^2-$part of the second fundamental form, i.e. $II=Q+Q^*$
for minimal surfaces. The straight line from $0$ to $1$ in $\CP^1$
corresponds to the geodesics $Q_1P_1,$ $Q_3P_1$ and $Q_5P_1$ in $S^3$ lying on $M.$ So the geodesic curvature of $Q_1P_1\subset S^3$ vanishes
which implies the assertion.
\end{proof}

\subsection{The holomorphic structures}
We use the symmetries to compute the spinor bundle $S\to M$ associated
to the Lawson genus $2$ minimal surface and the holomorphic structure $\nabla''$ on $V.$
\begin{Pro}\label{isospin}
Let $f\colon M\to S^3$ be a conformal minimal immersion. Let
$\Psi\colon S^3\to S^3$ and $\psi\colon M\to M$ be orientation
preserving isometries such that
$f\circ\psi=\Psi\circ f.$
Let $S\to M$ be the spinor bundle associated to $f.$
Then $$\psi^*S=S$$
as holomorphic bundles.
\end{Pro}
\begin{proof}
Because $S^3$ is simply-connected there is only
one $Spin-$structure on $S^3.$ Hence, 
$$(\Psi^*V,\Psi^*\nabla^{Spin})\cong(V,\nabla^{Spin}).$$
From the assertion one sees
$$\psi^*V=\psi^*f^*V=f^*\Psi^*V\cong f^*V=V$$
as unitary bundles with unitary connections on $M.$ Now $S^{-1}$ is 
the $-i$ eigenbundle of the complex quaternionic 
structure $\mathcal J$ induced by left multiplication with 
$-R.$ But $-R$ is invariant under $\psi.$ The holomorphic structure of $S^{-1}\subset V$ is given by $\nabla'',$ which is also invariant under $\psi.$ 
Therefore, the holomorphic structure of $S^{-1}$ is invariant under $\psi.$
\end{proof}
The correspondence between equivalence classes of divisors and holomorphic line bundles is classical, see \cite{GriHa}.
As above, we denote by $Q_1,..,Q_6$ the Weierstrass points of $M.$ Because $g(M)=2$
there are exactly
$\# H^1(M;\Z_2)=16$ different spin structures on $M,$ i.e. 
holomorphic line bundles $L$ satisfying $L^2=K.$
We list all of them below: On the left side are the different spin bundles, and on the right side are their 
pullbacks under symmetry $\Phi_3$ (for the computations we use that $2Q_i-2Q_j,$
$i,j=1,..,6,$
and $Q_1+Q_3+Q_5-Q_2-Q_4-Q_6$ are principal divisors):
\begin{equation}\label{spinbundlelist}
\begin{array}{|c|c|}
  \hline
  L & \Phi_3^*L \\ 
  \hline
  L(Q_1) & L(Q_5)  \\
  L(Q_2) & L(Q_6)  \\
  L(Q_3) & L(Q_1)  \\
  L(Q_4) & L(Q_2)  \\
  L(Q_5) & L(Q_3)  \\
  L(Q_6) & L(Q_4)  \\
  L(Q_2+Q_3-Q_1) & L(Q_6+Q_1-Q_5)=L(Q_5+Q_6-Q_1)  \\
  L(Q_2+Q_4-Q_1) & L(Q_6+Q_2-Q_5)=L(Q_3+Q_4-Q_1)  \\
  L(Q_2+Q_5-Q_1) & L(Q_6+Q_3-Q_5)=L(Q_2+Q_4-Q_1)  \\
  L(Q_2+Q_6-Q_1) & L(Q_6+Q_4-Q_5)=L(Q_2+Q_3-Q_1)  \\
  L(Q_3+Q_4-Q_1) & L(Q_1+Q_2-Q_5)=L(Q_2+Q_5-Q_1)  \\
  L(Q_3+Q_5-Q_1) & L(Q_1+Q_3-Q_5)=L(Q_3+Q_5-Q_1)  \\
  L(Q_3+Q_6-Q_1) & L(Q_1+Q_4-Q_5)=L(Q_4+Q_5-Q_1)  \\
  L(Q_4+Q_5-Q_1) & L(Q_2+Q_3-Q_5)=L(Q_4+Q_6-Q_1)  \\
  L(Q_4+Q_6-Q_1) & L(Q_2+Q_4-Q_5)=L(Q_3+Q_6-Q_1)  \\
  L(Q_5+Q_6-Q_1) & L(Q_3+Q_4-Q_5)=L(Q_2+Q_6-Q_1)  \\
  \hline
  \end{array}
\end{equation}
From this table one gets that the only $\Phi_3-$invariant spinor bundle
is $L(Q_1+Q_3-Q_5).$
Because $\Phi_3$ satisfies the conditions of proposition \ref{isospin} we obtain
\begin{The}\label{S=}
The spinor bundle $S\to M$ of the Lawson genus $2$ 
surface is given by $$S=L(Q_1+Q_3-Q_5).$$
\end{The}

The holomorphic structure $\nabla''$ is given by $\dbar=\dvector{\dbar^* & -\frac{i}{2}Q^*\\0 & \dbar}$ on $V=S^{-1}\oplus S,$ where $\dbar$ and $\dbar^*$ are the holomorphic structures on $S$ and $S^{-1}$ given by theorem \ref{S=}, and $Q^*\in\Gamma(M;\bar KK^{-1})$ is the adjoint of the Hopf differential. 
A holomorphic bundle over a Riemann surface of rank $2$ and degree $0$ is called stable if it does not contain proper holomorphic subbundles of degree greater or equal $0.$ We refer to \cite{NR} for details about extensions
and stable bundles. 
Because $Q^*\in\Gamma(M;\bar KK^{-1})$ is not in the image of the corresponding $\dbar-$operator,  one sees that there are no holomorphic subbundles
of positive degree.
By \cite{NR}, $V$ is non-stable if and only if there exits a point $x\in M$ such that $Q^*\otimes s_x\in\Gamma(M;\bar KK^{-1}L(x))$ is in the image of the corresponding 
$\dbar-$operator. Here $s_x\in H^0(L(x))$ is the canonical section of $L(x)$ which has exactly a simple zero at $x.$ By Serre duality, this condition is satisfied exactly in the case, that $$\int_M(Q^*\otimes s_x,\alpha)=0$$ for all $\alpha\in H^0(K^2L(-x)).$ Otherwise said, $Q^*$ is perpendicular to the $2-$dimensional subspace of holomorphic quadratic differentials which have a zero at some arbitrary but fixed point $x\in M$ if and only if $V$ is non-stable.
Let $P_1,..,P_4$ be the umbilics of the Lawson surface, and $\omega_1,\omega_2\in H^0(M;K)$
be the hyper-elliptic differentials with $(\omega_1)=P_1+P_3,$ $(\omega_2)=P_2+P_4.$ Using the hyperelliptic picture of $M$ and the symmetries of the Lawson surface one can easily compute
\[\int_M(Q^*,\omega_1^2)=\int_M(Q^*,\omega_2^2)=0,\] and
\[\int_M(Q^*,\omega_1\omega_2)\neq0.\]
Therefore, the space of holomorphic quadratic differentials which are perpendicular to $Q^*$ has no common zero. We have proven 
\begin{The}
The holomorphic rank $2$ bundle $(V,\nabla'')$ associated to the Lawson genus $2$ surface is stable.
\end{The}
\begin{Rem}
The holomorphic structure of a bundle $V$ in a short exact sequence $0\to S^{-1}\to V\to S$ is determined by the line of its extension class $[-\frac{i}{2}Q^*]\in PH^1(K^{-1}.)$ 
This line is already determined by $\int_M(Q^*,\omega_1^2)=\int_M(Q^*,\omega_2^2)=0$ and
$\int_M(Q^*,\omega_1\omega_2)\neq0.$
\end{Rem}
\begin{Rem}
This theorem shows that the method of \cite{He} to get a global DPW potential works for the Lawson surface. In order to get more informations 
about this potential, we will go another way in section \ref{DPWpotential}.
\end{Rem}

\section{The holonomy of the Lawson surface}\label{sec:holonomy}
We want to study the effect of the symmetries of the Lawson surface on the associated family of flat connections
$\nabla^\zeta$ and its holonomy representation. 
\subsection{A family of flat connections on the $4-$punctured sphere}
We start with the $\Z_3-$symmetry and consider
the threefold covering $\pi\colon M\to M/\Z_3=\CP^1.$
We show that the family of flat connections can be pulled back from the quotient.

There exists a square root $\sqrt{d\pi}$ of 
$d\pi\in H^0(M;\Hom(\pi^*K_{\CP^1}, K_{M})).$ 
To see this note that
$\pi^*K_{\CP^1}=L(-2Q_1-2Q_3-2Q_5),$ and that
$d\pi$ is determined by
$$s_{-2Q_1-2Q_3-2Q_5}\mapsto c s_{2P_1+..+2P_4-2Q_1-2Q_3-2Q_5}$$
for some nonzero constant $c.$ 
The pullback of the spinor bundle $\mathcal O(-1)\to\CP^1$
is given by $\pi^*\mathcal O(-1)=L(-Q_1-Q_3-Q_5).$
We define the holomorphic map 
$$\sqrt{d\pi}\colon\pi^*\mathcal O(-1)\to S=L(Q_1+Q_3-Q_5)$$
by the property
$$\sqrt{d\pi}(s_{-Q_1-Q_3-Q_5})=\sqrt cs_{P_1+..+P_4-Q_1-Q_3-Q_5}$$
for a fixed square root $\sqrt c$ of $c.$
This gives a commuting diagram
where the vertical maps are the spin double coverings:
\bigskip

$\xy
(80,0)*{}="A"; 
(102,0)*{}="B";
(80,-18)*{}="C";
(102, -18)*{}="D";
(11,-9)*{\,}; 
"A"*{\pi^*\mathcal{O}(-1)}; 
"B"*{S}; 
"C"*{\pi^*K_{\CP^1}};
"D"*{K_M};  
{\ar^{\sqrt{d\pi}} "A"+(9,0); "B"+(-2,0)}
{\ar_{d\pi} "C"+(7,0); "D"+(-4,0)}
{\ar_{} "A"+(-2,-3); "C"+(-2,3)}
{\ar^{} "B"+(-3,-3); "D"+(-3,3)}
\endxy
$
\bigskip

We also consider the inverse map 
$$\sqrt{d\pi}^{-1}\colon\pi^*\mathcal O(1)\to S^{-1}$$  
which is meromorphic with simple poles at $P_1,..,P_4$
and satisfies $\sqrt{d\pi}^{-1}\otimes\sqrt{d\pi}=1\in\C.$
Altogether we obtain a meromorphic map
\begin{equation}\label{Psidiag}
\Psi\colon \pi^*\mathcal O(1)\oplus\pi^*\mathcal O(-1)\to S^{-1}\oplus
S.
\end{equation}
The induced map on the determinant line bundles is an isomorphism.

Next, we consider the Hopf differential
$Q\in H^0(M;K^2)\subset H^0(M;K\End(S^{-1}\oplus S))$ and the meromorphic quadratic differential
$$\hat Q=\frac{ir}{z^4-1}(dz)^2\in \mathcal M(\CP^1;K_{\CP^1}^2)\subset 
\mathcal M(\CP^1;K_{\CP^1}\End(\mathcal O(1)\oplus\mathcal O(-1))),$$
where $r$ is as in proposition \ref{LawsonRS}.
Then its pull-back as an endomorphism-valued $1-$form 
$$\tilde Q:=\pi^*\hat Q\circ \pi\in \mathcal M(M;K\pi^*K_{\CP^1})\subset 
\mathcal(M;K\End(\pi^*\mathcal O(1)\oplus\pi^*\mathcal O(-1)))$$
has simple poles at $P_1,..,P_4.$
By construction, we have $$\Psi\tilde Q=Q\Psi.$$
Similarly, we can also define
 a Higgsfield on $\CP^1,$ which corresponds via pullback
and $\Psi$ to the Higgs field of the Lawson surface:
Consider
$$1=:\hat\Phi\in H^0(\CP^1;K_{\CP^1}\End(\mathcal O(1),\mathcal O(-1))).$$
and its pull-back
$$\tilde\Phi:=\pi^*\hat\Phi\in H^0
(M;K\End(\pi^*\mathcal O(1),\pi^*\mathcal O(-1)))$$
with double zeros at $P_1,..,P_4.$
Again, we have
$$\Psi\tilde\Phi=\Phi\Psi$$
on $M.$

As the Lawson surface has a $\Z_3-$symmetry the Riemannian metric on $M$ is invariant under $\Z_3$ and
induces a Riemannian metric $g$ on 
$\CP^1$ with conical singularities at the points $p_1,..,p_4,$
see \cite{Tr}:
The metric can be written
as a multiple of the constant curvature $1$ metric $g_0$ on $\CP^1,$ i.e.
$$g=e^{2\lambda}g_0$$ for some function 
$\lambda\colon\CP^1\setminus\{p_1,..,p_4\}\to\R.$
Because the pullback metric $\pi^*g$ is smooth on $M$ and $\pi$
has branch order $2$ at $P_i=\pi^{-1}(p_i)$ one sees that
$\lambda$ has a $\log$ singularity of order $-\frac{2}{3}.$ This means 
with respect to a holomorphic chart $z$ centered at $p_i$ the function
$\lambda$ can be written as
$$\lambda(z)=-\frac{2}{3}\ln{|z|}+f(z)$$ 
where $f$ is a locally defined smooth function.
The Riemannian metric $g$  on $\CP^1$ induces unitary metrics on 
$\mathcal O(\pm1)$ and on $\pi^*\mathcal O(\pm1)$ with conical singularities. The conformal factor
with respect to a smooth unitary metric has $\log-$singularity of order 
$\mp\frac{1}{3}$ on $\mathcal O(\pm1)$ and of order
$\mp1$ on $\pi^*\mathcal O(\pm1).$

With respect to 
these unitary metrics there exist
the adjoint operators
$$-\frac{i}{2}\hat Q^*\in\Gamma(\CP^1\setminus\{p_1,..,p_4\};
\bar K_{\CP^1}\End(\mathcal O(1)\oplus\mathcal O(-1)))$$ 
and
$$\hat\Phi^*\in\Gamma(\CP^1\setminus\{p_1,..,p_4\};
\bar K_{\CP^1}\End(\mathcal O(1)\oplus\mathcal O(-1)))$$
of $-\frac{i}{2}\hat Q$ and $\hat \Phi.$ They are $\bar K-$forms with values in $\End_0(\mathcal O(1)\oplus\mathcal O(-1)).$
Their pull-backs as endomorphism-valued $1-$forms are the adjoints of $-\frac{i}{2}\tilde Q$ and $\tilde \Phi$ by construction:
$$\pi^*(-\frac{i}{2}\hat Q^*)=\frac{i}{2}\tilde Q^*\in\Gamma(M\setminus\{P_1,..,P_4\};
\bar K_{M}\End(\pi^*\mathcal O(1)\oplus\pi^*\mathcal O(-1)))$$ 
and
$$\pi^*\hat\Phi^*=\tilde\Phi^*\in\Gamma(M\setminus\{P_1,..,P_4\};
\bar K_{M}\End(\pi^*\mathcal O(1)\oplus\pi^*\mathcal O(-1))).$$
They satisfy $\Psi\tilde Q^*=Q^*\Psi$ and  $\Psi\tilde\Phi^*=\Phi^*\Psi.$ 

What we have see is, that there is a holomorphic unitary bundle $\mathcal O(-1)\oplus\mathcal O(1)$ on the $4-$punctured sphere $\CP^1\setminus\{p_1,..,p_4\}$ together with a Higgs field, a holomorphic quadratic differential and their adjoint operators, which are mapped via the pullback $\pi^*$ and the meromorphic isomorphism $\Psi$ to the corresponding data of the Lawson surface. 
We need to define connections on
$\mathcal O(\pm1)\to\CP^1.$ A connection $\nabla$ with poles
(on a line bundle) is the sum of an ordinary connection with a meromorphic $1-$form $\omega.$ The residuum of $\nabla$ at a point
$p$ is given by $res_p\nabla:=res_p(\omega).$ 
\begin{Lem}\label{pushspincon}
There exists a unitary connection $\nabla^{\mathcal O}$ on 
$\mathcal O(-1)\to\CP^1$
with poles at $p_1,..,p_4$ such that the induced connection on $\mathcal O(-1)\otimes\mathcal O(-1)=K\to\CP^1\setminus\{p_1,..,p_4\}$ is the Levi-Civita
connection of the Riemannian metric with singularities at 
$p_1,..,p_4.$ The pullback $\pi^*\nabla^{\mathcal O}$ is a connection
on $\pi^*\mathcal O(-1)$ such that $\sqrt{d\pi}$ is parallel in 
$\Hom(\pi^*\mathcal O(-1)),S)$ where $S$ is equipped with the spin connection on $M.$ The residuums are given by $res_{p_i}\nabla=-\frac{1}{3}.$
 \end{Lem}
\begin{proof}
It is enough to show that the Levi-Civita connection $\nabla^K$ on $K_{\CP^1}$
with respect to the metric with singularities is a connection with poles
such that $res_{p_i}\nabla^K=-\frac{2}{3}.$ To see this we mention that complex connections on a line bundle $S$ and on the line bundle $S^2$
are in 1:1 correspondence. If one fixes corresponding connections on
$S$ and on $S^2,$ the connection forms for any other pair of
corresponding connections differ by the factor $2.$

The Levi-Civita connections of conformally equivalent metrics
$g=e^{2\lambda}g_{round}$ and $g_{round}$ differ on $K$ by
the the form $-2\del\lambda=-(d\lambda-i*d\lambda)\in\Gamma(K).$ 
The singularities of the metric on $\CP^1$ are given in such a way
that the conformal factor $\lambda$ has $\log-$singularities of order $-\frac{4}{3}.$
Thus $\del\lambda$ consists of the sum of a $C^\infty-$form and
a meromorphic $1-$form.
The residuums can be computed via the degree formula for
connections with poles by using the symmetries of the Riemannian metric.
\end{proof}
The induced connections on the dual bundle $\mathcal O(1)$ and on the
pullback bundles $\pi^*\mathcal O(\pm1)$ are also denoted by $\nabla^{\mathcal O}.$
We obtain
\begin{Pro}
There is an unitary connection $$\nabla^{\CP^1}:=\dvector{\nabla^{\mathcal O} & -\frac{i}{2}\hat Q\\
-\frac{i}{2}\hat Q^* & \nabla^{\mathcal O}}$$ on $\mathcal O(1)\oplus\mathcal O(-1)\to\CP^1\setminus\{p_1,..,p_4\}$ such that the map $\Psi$ is parallel
aside from $P_1,..,P_4\in M$
with respect to $\pi^*\nabla^{\CP^1}$ on $\pi^*\mathcal O(1)\oplus\pi^*\mathcal O(-1)\to M\setminus\{P_1,..,P_4\}$ and the hypersurface
connection $\nabla$ on $V\to M.$ 
\end{Pro}
Consider the holomorphic family
$$\zeta\in\C\setminus\{0\}\mapsto\nabla^\zeta:= \nabla^{\CP^1}+\zeta^{-1}\hat\Phi-\zeta\hat\Phi^*$$ 
of connections  on $\mathcal O(1)\oplus\mathcal O(-1)\to\CP^1\setminus\{p_1,..,p_4\}$
with singularities at $p_1,..,p_4.$  
 Because of the construction the 
map $\Psi$ is parallel with respect to the pullback connections
$\pi^*\nabla^\zeta$
on $\pi^*(\mathcal O(1)\oplus\mathcal O(-1))$ and the corresponding flat connections $\nabla^\zeta$ 
given by equation \ref{nablafamily} on $V$ for all $\zeta\in\C^*.$ As $\pi$ is a holomorphic covering this implies that the connections $\nabla^\zeta$ on $\mathcal O(1)\oplus\mathcal O(-1)\to\CP^1\setminus\{p_1,..,p_4\}$ are flat for all $\zeta\in\C^*.$ Consequently, for any $z\in\CP^1\setminus\{p_1,..,p_4\},$ one obtains a holonomy representation
$$\mathcal H_z\colon\C^*\times\pi_1(\CP^1\setminus\{p_1,..,p_4\},p)\to\SL((\mathcal O(1)\oplus\mathcal O(-1))_z).$$
Because $z=0$ is a fix point for most of the symmetries, we consider the holonomy at $z=0.$ The first fundamental group of $\CP^1\setminus\{p_1,..,p_4\}$
is generated by the closed loops 
$\gamma_k\colon S^1\to\CP^1\setminus\{p_1,..,p_4\}$ centered
at $0\in\CP^1$ which go counter-clockwise around the points $p_k.$ They satisfy $\gamma_1*..*\gamma_4=1\in\pi_1(\CP^1\setminus\{p_1,..,p_4\},0).$
We fix a special unitary basis $e_1,e_2$ of $(\mathcal O(1)\oplus\mathcal O(-1))_{|0}$ and corresponding bases of $V_{|Q_i},$ $i=1,3,5.$
Let 
$$\zeta\mapsto H_k(\zeta)\in\SL(2;\C)$$
be the holonomy of $\nabla^\zeta$ around $\gamma_k$ centered
at $0$ with respect to the basis. 
The lift  of $\gamma_k^3$ is a trivial loop on $M$ for $k=1,..,4.$
Therefore $$H_k^3(\zeta)=\Id$$ as the
connection $\nabla^\zeta$ with singularities at $p_1,..,p_4$ can be desingularised on
$V\to M.$ 
The conjugacy class of  $H_k$ is independent of $\zeta$ for $k=1,..,4.$
If $H_k$ would be the identity, the symmetries of the Lawson surface
would imply that the holonomy of the Lawson surface is trivial for all $\zeta,$ see theorem \ref{holonomyupdown} and propositions \ref{Hpm2} and \ref{Hpm1} below. But this can only happen for the round sphere,
see \cite{H1}. Therefore, the eigenvalues of $H_k$ are $\xi:=e^{\frac{2}{3}\pi i},\,\xi^2\in\C.$ This is equivalent to
$$\Tr(H_k)=\xi+\xi^2=-1,$$ independently of $\zeta.$

Consider 
the following closed oriented geodesic polygons on $M$ centered at $Q_1:$
\begin{equation}\label{canonicalbasis}
\begin{split}
\Gamma_1:=&Q_1P_1Q_3P_2Q_1\\
\Gamma_2:=&Q_1P_1Q_3P_4Q_1\\
\Gamma_3:=&Q_1P_1Q_5P_4Q_1\\
\Gamma_4:=&Q_1P_1Q_5P_2Q_1.\\
\end{split}
\end{equation}
They generate the first fundamental group of $M.$ 
In fact $a_1:=\Gamma_1*\Gamma_4^{-1}, b_1:=\Gamma_2, a_2:=\Gamma_3, b_2:=\Gamma_2$ is a canonical basis of $\pi_1(M).$

We denote
by $$\zeta\mapsto A_k(\zeta)\in\SL(2,\C)$$ 
the holonomy of $\nabla^\zeta$ on $V$ around the $\Gamma_k$
centered at $Q_1$ with respect to the basis of $V_{|Q_1}$ chosen above. 
\begin{The}\label{holonomyupdown}
There is a family of flat connections $\nabla^\zeta$ with singularities at 
$p_1,..,p_4$
on $W\to\CP^1$ which correspond under pull-back and the map $\Psi$ to the family of flat connections $\nabla^\zeta$ associated to Lawson's genus $2$ surface.
The holonomies are related by the following formula:
 \begin{equation}
\begin{split}
A_1=&H_2^{-2}\circ H_1=H_2\circ H_1\\
A_2=&H_4^{-2}\circ H_1=H_4\circ H_1\\
A_3=&H_4^{-1}\circ H_1^2=H_4^2\circ H_1^2\\
A_4=&H_2^{-1}\circ H_1^2=H_2^2\circ H_1^2.\\
\end{split}
\end{equation}
\end{The}
\begin{proof}
We have already proven the first part. Therefore, the assertion regarding the holonomies follows from the following observation:
A lift of the loop $\gamma_1\colon S^1\to \CP^1\setminus\{p_1,..,p_4\}$ to $M$ with starting point
$Q_1$ is a curve with end point $Q_3$ which is homotopic
(with fixed endpoints) to the geodesic polygon
$Q_1P_1Q_3,$ and a lift of the loop $\gamma_2\colon S^1\to\CP^1\setminus\{p_1,..,p_4\}$ with starting point
$Q_1$ corresponds to a curve with end point $Q_5$ which is homotopic
(with fixed endpoints) to the geodesic polygon
$Q_1P_2Q_5.$ 
Analogous statements hold for 
$\gamma_3,\gamma_4.$ 
\end{proof}

\subsection{Symmetries and Holonomy}\label{sec:symmetries_holonomy}

We first 
consider the orientation preserving isometry
$\Phi_2\colon S^3\to S^3; (a,b)\mapsto(a,-b).$ This induces the hyperelliptic involution of the Lawson surface 
 with fixpoints
$Q_1,..,Q_6.$ The corresponding map on the $3-$fold covered $\CP^1$ is $z\mapsto-z.$
Then 
\begin{equation}\label{fphi2}
f\circ\Phi_2=\Phi_2\circ f=\dvector{- i & 0\\ 0 & i}f\dvector{ i & 0\\ 0 & -i},
\end{equation}
where the last equation is in $S^3=\SU(2,\C)\subset\SL(2,\C).$
The connections $\nabla^\zeta$ on $V\to M$ and $\nabla^\zeta$ on $\mathcal O(1)\oplus\mathcal O(-1)\to\CP^1\setminus\{p_1,..,p_4\}$ are invariant
under $\Phi_2\colon M\to M$ and $z\mapsto-z$ on $\CP^1,$ respectively. Due to equation \ref{fphi2}, the chosen basis 
of $(\mathcal O(1)\oplus\mathcal O(-1))_{|0}$ changes, and we obtain

\begin{Pro}\label{Hpm2}
The holonomies around $p_k$ and $p_{k\pm2}$ satisfy
\begin{equation*}
\begin{split}
H_3(\zeta)&=\dvector{-i & 0\\0 &i}H_1(\zeta)\dvector{i & 0\\0 &-i}\\
H_4(\zeta)&=\dvector{-i & 0\\0 &i}H_2(\zeta)\dvector{i & 0\\0 &-i}\\
\end{split}
\end{equation*}
for all $\zeta\in \C^*.$
\end{Pro}
Next we consider the isometry $\tau:=\gamma_{S_2}\circ\gamma_{P_1Q_1},$ which is holomorphic on the surface. On $\CP^1$ the induced action is given by $z\mapsto iz.$ But this symmetry changes orientation in space.
This implies that the Hopf field changes by the factor $-1,$ i.e. $\tau^*Q=-Q.$ Of course, the spin connections on the diagonal of $\nabla^\zeta$ and
the Higgs field are invariant under $\tau,$ like their adjoint operators. Therefore, we obtain
\begin{equation}\label{nablaHpm1}
\begin{split}
\tau^*\nabla^\zeta=\nabla^{-\zeta}\cdot\dvector{ i & 0\\ 0 & -i},
\end{split}
\end{equation}
and $H_1(-\zeta)$ and $H_2(\zeta)$ must be related:
\begin{Pro}\label{Hpm1}
The holonomies $H_1$ and $H_2$ satisfy
\begin{equation*}
\begin{split}
H_2(\zeta)=\dvector{e^{-i\frac{\pi}{4}} & 0 \\ 0 & e^{i\frac{\pi}{4}}}H_1(-\zeta)\dvector{e^{i\frac{\pi}{4}} & 0 \\ 0 & e^{-i\frac{\pi}{4}}}
\end{split}
\end{equation*}
for all $\zeta\in\C^*.$
\end{Pro}
\begin{proof}
The surface obtained from the family of connections $\zeta\mapsto\nabla^{-\zeta}$ with respect to the same basis of 
$V_{|Q_0}\cong\C^2$ as for $f,$ is exactly
the surface $f^{-1}\colon  M\to S^3.$
But 
$$\dvector{e^{i\frac{\pi}{4}} & 0 \\ 0 & e^{-i\frac{\pi}{4}}}\tau(f)\dvector{e^{-i\frac{\pi}{4}} & 0 \\ 0 & e^{i\frac{\pi}{4}}}=f^{-1}.$$
and with equation \ref{nablaHpm1} we obtain the assertion.
\end{proof}
Next we deal with the orientation preserving symmetry $\gamma_{P_1Q_1}\colon S^3\to S^3, (a,b)\mapsto(\bar a,\bar b),$ which is 
antiholomorphic on the surface. The induced action on $\CP^1$ is $\gamma_{P_1Q_1}(z)=\bar z.$
\begin{Pro}\label{zbar-symmetry}
The holonomy around $P_1$ satisfies
\begin{equation*}
\begin{split}
H_1^2(\zeta^{-1})=H_1^{-1}(\zeta^{-1})=\dvector{0 & -1\\ 1& 0}H_1(\zeta)\dvector{0 & 1\\ -1& 0}.
\end{split}
\end{equation*}
for all $\zeta\in\C^*.$
\end{Pro}
\begin{proof}
We have $$\tilde f:=f\circ\gamma_{P_1Q_1}=\dvector{0 & -1\\ 1& 0}f\dvector{0 & 1\\ -1& 0}$$ as unoriented maps. The oriented normal of $\tilde f$ is $-N\circ\gamma_{P_1Q_1},$ where $N\colon M\to S^2$ is the normal of $f.$  Then, $\tilde f$ with oriented normal 
$-N\circ\gamma_{P_1Q_1}$ has holonomy 
$$\tilde H_1(\zeta)=\dvector{0 & -1\\ 1& 0}H_1(\zeta)\dvector{0 & 1\\ -1& 0}.$$
We need to figure out which holonomy is given by $\tilde f$ with respect 
to the antiholomorphic orientation:
To do so we consider the associated family of connections.
Note that $\gamma_{P_1Q_1}^*K=\bar K$ and 
$\gamma_{P_1Q_1}^* S=\bar S=S^{-1},$ where we use the unitary metric for identification.
We denote the Maurer-Cartan form of $\tilde f$ by $\tilde \phi,$ and so on.
 Then the $0th$-order part of the families of connections $\nabla^\zeta$ and $\tilde\nabla^\zeta$ are the pull-backs under $f$ and $f\circ\gamma_{P_1Q_1},$ respectively,
of the Riemannian spin connection on the spinor bundle over $S^3.$ Therefore, the $0th$-order part of the family of connections of $f$ pulls back to the
$0th$-order part of the family of connections of $f\circ\gamma_{P_1Q_1}.$ Moreover $$\gamma_{P_1Q_1}^*\phi=(f\circ\gamma_{P_1Q_1})^{-1}d(f\circ\gamma_{P_1Q_1}).$$
and we get 
$$\gamma_{P_1Q_1}^*(\Phi-\Phi^*)=\tilde\Phi-\tilde\Phi^*$$
as $1-$forms with values in $\End_0(V).$ 
We obtain
$$\gamma_{P_1Q_1}^*\Phi=-\tilde\Phi^*\,\,\,\,\,\, \text{and}\,\,\,\,\,\, \gamma_{P_1Q_1}^*\Phi^*=-\tilde\Phi.$$
This implies that
$$\gamma_{P_1Q_1}^*\nabla^{\zeta^{-1}}$$ is
the  family of connections associated to $\tilde f=f\circ\gamma_{P_1Q_1}.$ Note that
$f\circ \gamma_{P_1Q_1}$ is the gauge between $\gamma_{P_1Q_1}^*\nabla^1$ and $\gamma_{P_1Q_1}^*\nabla^{-1}$ with respect to the same basis of $V_{|Q_1}$ as for $f.$
Moreover the symmetry $\gamma_{P_1Q_1}(z)=\bar z$ maps the curve $\gamma_1$
which goes around $1\in\C$ counterclockwise to a curve homotopic to $\gamma_1^{-1}$ which goes around $1\in\C$ clockwise.
This implies the statement.
\end{proof}

There exists another relation between the holonomies: The product
$$H_4H_3H_2H_1=\Id$$ is the identity, as the family of connections is well-defined on the $4-$punctured sphere.
Then Proposition \ref{Hpm2} yields
 $$\dvector{-i & 0\\ 0 & i}H_2H_1\dvector{i & 0\\ 0 & -i} =(H_2H_1)^{-1}.$$

\section{A DPW potential for Lawson's genus 2 surface}\label{DPWpotential}
The idea of the DPW method is to gauge the family $\nabla^\zeta$
(see equation \ref{nablafamily}) into a family of meromorphic connections in a way which can be
reversed. In principle, one can construct all minimal surfaces in $S^3$ by this method. But in concrete situations, it is very difficult to produce surfaces with prescribed properties. For example, there is no compact minimal surface of genus $g\geq2$ constructed via the DPW method
up to now. Nevertheless, there is some work of the author \cite{He} which shows, that the DPW method should work fine for compact
surfaces of genus $2:$ There it is proven that there exists a globally defined DPW potential which gives back the minimal surface. 
Here, we consider the special case of the Lawson genus $2$ surface,
and we can show (theorem \ref{explicitform}) the existence of a globally defined DPW potential, whose behavior on the surface is completely described. The freedom of the potential is given by two unknown functions in $\zeta.$  

\begin{Def}
A meromorphic connection $\nabla$ on a holomorphic vector bundle $(V,\dbar)$ is a connection with 
singularities which can be written with respect to a local holomorphic frame as $d+w$ where $w$ is an meromorphic endomorphism-valued $1-$form.
\end{Def}

Of course, on Riemann surfaces meromorphic connections are flat. 
For line bundles
there exists the degree formula
$$res(\nabla)=-deg(L)$$ on Riemann surfaces, where the $res(\nabla)$ is the sum of all local 
residui.
\begin{The}\label{DPWex}
Let $\nabla^\zeta$ be the holomorphic family of flat connections on $V$ associated to Lawson's genus 2 surface $f\colon M\to S^3.$ 
Let $Q_1,..,Q_6$ be the Weierstrass points of $M$. 
Then there exists a holomorphic map
$$\tilde B:\zeta\in B(0;\epsilon)\subset\C\to \Gamma (M\setminus\{Q_1,..,Q_6\},\End(V))$$ which
satisfies $B_0=\dvector{1 & * \\ 0 & 1}$ and $\det B_\zeta=1$ for all $\zeta,$ such that the gauged connection
$$\nabla^\zeta\cdot B_\zeta$$ is a holomorphic family of meromorphic connections $\hat\nabla^\zeta$ for 
$\zeta\in B(0;\epsilon)\setminus\{0\}\subset\C$ on the (fixed) holomorphic vector bundle $(V=S^{-1}\oplus S,\dbar^{spin}).$

More precisely the family has an expansion
$$\hat\nabla^\zeta=\dvector{\nabla_{0}^* & \zeta^{-1}\\ -\frac{i}{2}Q & \nabla_{0}}+higher\ order\ terms,$$
for some meromorphic connection $\nabla_0$ on $S$ and the Hopf field $Q\in H^0(K^2)$ of the surface.
The connections have poles of order $1$ on the diagonal at $Q_1,..,Q_6$ and of order $2$ in the upper right and lower left corner at $Q_2,Q_4,Q_6$ respectively 
$Q_1,Q_3,Q_5.$
\end{The}
\begin{proof}
The condition that $\nabla^\zeta\cdot B_\zeta$ is a holomorphic family of meromorphic connections
on the holomorphic bundle $S^{-1}\oplus S$ translates easily to
\begin{equation}\label{dbarB}
\dbar^{spin}B=(\frac{i}{2}Q^*+\zeta\Phi^*) B,
\end{equation}
with $Q^*\in\Gamma(\bar KK^{-1})$ and $\Phi^*\in\Gamma(\bar KK).$
Writing $$B=\sum_{k\geq0} B_k\zeta^k$$ and
$$B_k=\dvector{ a_k & b_k\\ c_k & d_k}$$ one obtains the equations
\begin{equation}
\begin{split}
\dbar a_k&=\frac{i}{2}Q^*c_k\\
\dbar c_{k+1}&=\Phi^*a_k\\
\dbar b_k&=\frac{i}{2}Q^*d_k\\
\dbar d_{k+1}&=\Phi^*b_k\\
\end{split}
\end{equation}
for $a_k,d_k\in\Gamma(M;\underline\C), b_k\in\Gamma(M;K^{-1})$ and $c_k\in\Gamma(M;K).$
If we would take $a_0=d_0=1,$ $c_0=0,$ then, by Serre duality, there does not
exists a smooth $b_0$ satisfying the equation above. To overcome this problem, we search for solutions with singularities.
Set
$$a_0=1,\, d_0=1,\, c_0=0.$$
Take the divisors $D=Q_1+Q_3+Q_5$ and
$\tilde D=Q_2+Q_4+Q_6$ which are invariant under the $\Z_2$ and under the $\Z_3$ action. Note that $L(D)=L(\tilde D)=KS.$ Now consider\begin{equation}
\begin{split}
\tilde a_k&=a_k\otimes s_D\in\Gamma(KS)\\
\tilde b_k&=b_k\otimes s_{\tilde D}\in\Gamma(S)\\
\tilde c_k&=c_k\otimes s_D\in\Gamma(K^2S)\\ 
\tilde d_k&=d_k\otimes s_{\tilde D}\in\Gamma(KS).
\end{split}
\end{equation}
 We get new equations
\begin{equation}
\begin{split}
\dbar\tilde a_k&=\frac{i}{2}Q^*\tilde c_k\in\Gamma(\bar KKS)\\
\dbar\tilde c_{k+1}&=\Phi^*\tilde a_k\in\Gamma(\bar KK^2S)\\
\dbar\tilde b_k&=\frac{i}{2}Q^*\tilde d_k\in\Gamma(\bar KS)\\
\dbar\tilde d_{k+1}&=\Phi^*\tilde b_k\in\Gamma(\bar KKS).\\
\end{split}
\end{equation}
Again, Serre duality tells us that there does always exist a solution for each of these equations. But we need more: we want the $\zeta-$series $\sum\tilde B_k\zeta^k$ to be convergent and $\det B_{\zeta}=1.$ We explain how this can be achieved. Note that all occurring bundles inherit canonical unitary metrics from the surface metric. These give us fixed Sobolev norms and spaces.
By Poincare inequality there exists a constant $c>0$ such that the solution $s$ for any of the above equations, which is unique by the property of being orthogonal to the kernel of the corresponding 
$\dbar-$operator, satisfies
$$\parallel s\parallel\leq c\parallel \dbar s\parallel.$$
Note that, if the right hand side of any of these equations is symmetric with
respect to the $\Z_2$ or $\Z_3$ symmetry, the unique solution has also this symmetry. 
We take always this unique solution, and solve for $\tilde a_{k+1},..,\tilde d_{k+1}$ inductively. Thus we obtain $\parallel\tilde B_k\parallel<C^k$ for some constant $C,$
which implies smooth convergence for small $\zeta.$
Set $$B:=\sum_{k\geq0} \dvector{\tilde a_k\otimes s_{-D} & 
\tilde b_k\otimes s_{-\tilde D} \\ \tilde c_k\otimes s_{-D} & \tilde d_k\otimes s_{-\tilde D} }\zeta^k.$$
Then $B$ satisfies equation \ref{dbarB}.
Because of this equation the determinant $\det(B_\zeta)$ is a meromorphic function on $M$ for all $\zeta$. Moreover it is invariant under the $\Z_2$ action.
 Clearly, $det(B_\zeta)\geq -D-\tilde D$ for all $\zeta.$ But the only $\Z_2-$invariant functions with this divisor inequality are the constants. Thus there exists $h\colon B(0;\epsilon)\to\C$ with $\det(B_\zeta)=h(\zeta).$ Therefore
 $$B_\zeta\dvector{\frac{1}{h(\zeta)} & 0\\ 0 & 1}$$
 is the gauge we were looking for.
\end{proof}
As we have seen in the proof of the previous theorem, the meromorphic connections $\hat\nabla^\zeta$ have some of the symmetries of the Lawson surface by construction. We use these symmetries to write down the corresponding DPW potential almost explicitly.
To do so, we trivialize $S^*\oplus S\to M\setminus\{Q_2,Q_4,Q_6,P_1,..,P_4\}$ using the meromorphic sections
 $s=s_{Q_2+Q_4+Q_6-P_1-P_2-P_3-P_4}\in\mathcal M(S^*)$ and 
 $t=s_{-Q_2-Q_4-Q_6+P_1+P_2+P_3+P_4}\in\mathcal M(S).$
\begin{The}\label{explicitform}
With respect to the meromorphic frame $(s,t)$ of $S^*\oplus S$
and up to a diagonal gauge only depending on $\zeta$
the family of connections given by theorem \ref{DPWex} can be written as
$d+\xi$ with
$$\xi=\pi^*\dvector{-\frac{4}{3}\frac{z^3}{z^4-1}+\frac{A}{z}  & \zeta^{-1}+Bz^2 \\ \frac{G}{(z^4-1)} +\frac{\zeta H}{z^2(z^4-1)}& \frac{4}{3}\frac{z^3}{z^4-1}-\frac{A}{z}}dz $$
for some $\zeta$ depending even functions $A,B,G,H$ which satisfy
$H= A+A^2$ and $B= -\frac{1}{G}(-\frac{1}{3}+A+(\frac{1}{3}-A)^2).$
\end{The} 
\begin{proof}
Note that the upper right corner with respect to the holomorphic decomposition $S^*\oplus S$ of $\nabla^\zeta$ has the invariant meaning of a meromorphic function on 
$M$ with at most double poles at $Q_2,Q_4,Q_6.$ There is a well-defined holomorphic function $h(\zeta)$ on $B(0;\epsilon)\setminus\{0\}$ which is the constant part (but depending on $\zeta$) of the upper right corner. From the starting condition $a_0=1=d_0$ and $c_0=0$ we get the Laurent expansion $h(\zeta)=\zeta^{-1}+h_0+...$. For $\epsilon$ small enough we can take a square root $g$ of $\zeta h(\zeta).$ Instead of working with $\hat\nabla^\zeta$ we gauge it by $\dvector{g(\zeta) &0\\0&1/g(\zeta)},$ so that the part of the upper right corner, which is constant along $M,$  is given by $\zeta^{-1}.$

With respect to the given trivialization the connection $1-$ form $\xi,$ the so-called DPW potential, is a $\mathfrak{sl}(2;\C)-$valued family of meromorphic $1-$forms. We want to deduce the symmetries of the potential $\xi$ from the symmetries of
the family of connections $\nabla^\zeta$ and of the gauge $B_\zeta.$ 
We start with the $\Z_3$ symmetries. Note that $\Phi_3^*(s,t)=(s,t).$ By construction the family of connections $\hat\nabla^\zeta$ is also invariant under $\Z_3.$ Thus $\xi$ is invariant under the $\Z_3-$action, i.e. $\Phi_3^*\xi=\xi.$ Similarly, the family of connections $\hat\nabla^\zeta$ is also invariant under $\Z_2.$
But as $\Phi_2^*(s,t)=( -i s, i t)$ the generator $\Phi_2$ of the $\Z_2-$action satisfies $$\Phi_2^*\xi=\xi\cdot\dvector{ i & 0 \\ 0 & -i}.$$

Next, we look at $\tau=\gamma_{S_2}\circ\gamma_{P_1Q_1}$ which is holomorphic on the surface but changes orientation in space. Note that $\tau^*(s,t)=( e^{-\frac{\pi i}{4}} s, e^{\frac{\pi i}{4}} t),$ and $\tau^*\nabla^\zeta=\nabla^{-\zeta}\dvector{i & 0 \\ 0 & -i}.$ From $\tau^*Q^*=-Q^*, \tau^*\Phi^*=\Phi^*$ and from the recursive construction of
$B$ one sees that
$$\tau^*B(\zeta) =\dvector{-i & 0 \\ 0 & i}B(-\zeta)\dvector{i & 0 \\ 0 & -i}.$$
Altogether we obtain the following symmetry
$$\tau^*\xi(\zeta)=\xi(-\zeta)\cdot \dvector{ e^{i\frac{\pi}{4}}  & 0 \\ 0 & e^{-i\frac{\pi}{4}} }$$
for the DPW potential $\xi.$

The symmetry $\tilde\tau=\gamma_{S_1}\circ\gamma_{P_1Q_1}$ is 
more difficult to handle. The reason for this is that it interchanges the 
points $Q_1,Q_3,Q_5$ and $Q_2,Q_4,Q_6.$ 
As it is holomorphic on the surface and orientation reversing in space, we have again $\tilde\tau^*\nabla^\zeta=\nabla^{-\zeta}\dvector{i & 0 \\ 0 & -i}.$ Moreover $\tilde\tau^*(s,t)=(-iz s,i \frac{1}{z}t).$ We claim that
the symmetry of the potential takes the following form:
$$\tilde\tau^*\xi(\zeta)=\xi(-\zeta)\cdot g\cdot \dvector{ z \circ\pi & 0 \\ 0 & \frac{1}{z\circ\pi} },$$
for some (globally defined) automorphism $g.$
In fact, $g$ is given by
$$g=\dvector{a_0+a_1z^2 & \frac{1}{z} b_1(z^4-1)\\ \frac{c_{-1}}{z} & \frac{d_{-1}}{z^2}+d_0}\circ\pi$$
for some $\zeta$ depending functions $a_0,..,d_0.$ This can be deduced from the fact that $\tilde\tau^*\xi(\zeta)$ can be obtained by the same method as $\xi,$ see the proof of theorem \ref{DPWex}, with the difference that the gauge has singularities in the first column at $Q_2,Q_4,Q_6$ and in the second column at $Q_1,Q_3,Q_5.$
We will not give the details as it turns out below that the symmetry $\tilde\tau$
does not give any new information for the potential $\xi.$

Next, we list the poles of the potential $\xi.$ It has 
\begin{itemize}
\item poles on the diagonal of order $1$ at $Q_1,..,Q_6;$
\item poles on the diagonal of order $1$ at $P_1,..,P_4$ with residuum $\mp1;$
\item poles in the lower left corner up to order $2$ at $Q_1,Q_3,Q_5;$
\item poles in the lower left corner up to order $2$ at $P_1,..,P_4;$
\item poles in the upper right corner up to order $4$ at $Q_2,Q_4,Q_6.$
\end{itemize}
Note that the poles at $P_1,..,P_4$ and the poles of order $4$ 
instead of $2$ in the upper right corner come from the chosen 
trivialization $(s,t).$ 

We have enough informations to determine the potential:
Because $\xi$ and the trivialization have both the $\Z^3-$symmetry,
$\xi$ is the pullback of an $\zeta-$depending $\mathfrak{sl}(2;\C)-$valued meromorphic $1-$form $\tilde\xi$ on $\CP^1.$ This $1-$form on $\CP^1$
can be considered as the connection $1-$form with respect to the frame $(s_{\infty},s_{-\infty})$ 
of a holomorphic family of connections on $\mathcal O(1)\oplus\mathcal O(-1)\to\CP^1,$ compare with section \ref{sec:holonomy}.
Clearly $\tilde\xi$ has corresponding symmetries and pole behavior, for example, the poles on the diagonal have residuum $\pm\frac{1}{3}$ at the $4$th order roots of $1.$

The most general form of a potential on $\CP^1$ with the symmetries
for $\Phi_2,\Phi_3$ and $\tau$ and singularities as above is up to a diagonal gauge only depending on $\zeta$
\begin{equation}\label{potentialons2}
\dvector{-\frac{4}{3}\frac{z^3}{z^4-1}+\frac{A}{z}  & \zeta^{-1}+Bz^2 \\ \frac{G}{(z^4-1)} +\frac{\zeta H}{z^2(z^4-1)}& \frac{4}{3}\frac{z^3}{z^4-1}-\frac{A}{z}}dz 
\end{equation}
for some $\zeta$ depending even functions $A,B,G,H.$ 

The singularities at $z=0$ and $z=\infty$ are apparent.
In fact, the construction of the potential shows that there must be a $\zeta-$depending meromorphic gauge $g$ which is diagonal at $\zeta=0,$ and lower triangular respectively upper triangular for general $\zeta,$ such that $g$ gauges the singularities at $z=0$ respectively $z=\infty$ away.
Then, a simple computation shows 
\begin{equation}\label{ABGHeq}
\begin{split}
H&= A+A^2\\
B&= -\frac{1}{G}(-\frac{1}{3}+A+(\frac{1}{3}-A)^2)
\end{split}
\end{equation}
for all $\zeta.$

Again, a short computation shows, that each potential of the form 
\ref{potentialons2} with functions $A,B,G,H$ satisfying the equations
\ref{ABGHeq} posses the symmetry for $\tilde\tau.$
\end{proof}
The next task would be to determine the functions $A$ and $G.$ As we have seen, they cannot be computed out of the symmetries. These functions satisfy more complicated equations. Namely, the DPW potential must be unitarizable. This means, that the holonomy representation of the family of connections
for some $\zeta-$depending starting condition must extend to $\C\setminus\{0\},$ and must be unitary
for $\zeta\in S^1\subset\C.$ This starting condition is called dressing.
If one gauges the singularity at $z=0$ away with an lower triangular gauge and then
starts integrating at  $z=0$, i.e. finding a parallel frame,  then the dressing must be diagonal as one can see from the symmetries. The holonomy depends transcendentally on $A$ and $G,$
so it is a very hard problem to determine the exact form of $A$ and $G.$ The space of representations of $\pi_1(\CP^1\setminus{p_1,..,p_4}$ in $\SL(2,\C)$ modulo conjugation with fixed conjugacy classes around $p_1,..,p_4$ is a cubic surface in $\C^3,$ see \cite{BG}. One can show, that for fixed $\zeta\in\C^*$ the holonomy representation depends on $A$ and $G$ independently. Thus, one obtains an open nonempty subset of all possible representations from the potential $\xi.$
But, one cannot obtain all possible representations as the proof of the 
following theorem shows.
\begin{The}
The family of meromorphic connections $\hat\nabla$ and the DPW potential $\xi$ do not extend to
the whole unit circle.
\end{The} 
\begin{proof}
As the gauge $B$ is well-defined on the surface $M$ and not multi-valued, the monodromy representations of $\nabla^\zeta$ and of $\hat\nabla^\zeta$ are equivalent for all $\zeta\in \C^*.$
If $\hat\nabla^{\pm 1}$ would exist, the monodromy would be trivial.
This means there would exist two linear independent meromorphic sections $v,w\in\mathcal M(M,S^*\oplus S)$ parallel with respect to $\hat\nabla^1.$ As $\hat\nabla^1$ has its only poles at $Q_1,..,Q_6$ the same holds for $v$ and $w.$
Write $v=x\oplus y$ with respect to $S^*\oplus S.$ From the special form of $\xi$ one sees that $x$ has simple poles at $Q_2,Q_4,Q_6$  and $y$ has simple poles at $Q_1,Q_3,Q_5.$
Therefore, $y$ is a constant multiple of the meromorphic section
$s_{-Q_1-Q_3-Q_5+P_1+..+P_4}=\frac{1}{z\circ\pi}t\in\mathcal M(M,S).$
The same argument holds for a decomposition of $w,$ which shows that a parallel frame $v,w$ would not be linear independent at the points
$P_1,..,P_4,$ which is a contradiction. 
\end{proof}


\end{document}